\documentclass[9pt,final,letterpaper,technote]{IEEEtran}
\usepackage{amsmath}
\usepackage{amssymb}

\numberwithin{equation}{section}

\def\theequation{\arabic{section}.\arabic{equation}}

\newtheorem{prop}{Proposition}
\newtheorem{asmp}{Assumption}

\newtheorem{defn}{Definition}

\newtheorem{exmp}{Example}

\newtheorem{rem}{Remark}

\begin{document}

\title{Sufficient Lie Algebraic Conditions for Sampled-Data Feedback Stabilizability of Affine in the Control Nonlinear Systems}

\author{J.~Tsinias and D.~Theodosis %
\thanks{Authors are with the Department of Mathematics, National Technical University of Athens, Zografou Campus 15780, Athens, Greece, email: jtsin@central.ntua.gr (corresponding author), dtheodp@central.ntua.gr. }} %
\maketitle

\begin{abstract}
For general nonlinear autonomous systems, a Lyapunov characterization for the possibility of semi-global asymptotic stabilizability by means of a time-varying sampled-data feedback is established. We exploit this result in order to derive a Lie algebraic sufficient condition for sampled-data feedback semi-global  stabilizability of affine in the control nonlinear systems with non-zero drift terms. The corresponding proposition constitutes an extension of the ``Artstein-Sontag'' theorem on feedback stabilization.
\end{abstract}


\section{Introduction}
\label{S:1}

Many significant results towards stabilizability of nonlinear systems by means of sampled-data feedback control have appeared in the literature (see for instance \cite{art:1}, \cite{art:2}, \cite{art:4}-\cite{art:12} and relative references therein). In the recent works \cite{art:17} and \cite{art:19}, the concept of \textit{Weak Global Asymptotic Stabilizability by Sampled-Data Feedback} (SDF-WGAS) is presented for systems:
\begin{equation}\label{1.1}
\begin{array}{c}
\dot{x}=f(x,u),\, \, (x,u)\in \mathbb{R}^{n}\times \mathbb{R}^{m},\\
f(0,0)=0 
\end{array}
\end{equation}
and various Lyapunov-like sufficient characterizations of this property are examined. Particularly,  in \cite[Proposition 2]{art:19}, a Lie algebraic sufficient condition for SDF-WGAS is established for the case of affine in the control systems 
\begin{equation}\label{1.2}
\begin{array}{c}
 \dot{x}=f(x)+ug(x),\; \; (x,u)\in  \mathbb{R}^{n}\times\mathbb{R},\\
   f\left(0\right)=0  
   \end{array}
\end{equation} 
This condition constitutes an extension of the well-known ``Artstein-Sontag'' sufficient condition for asymptotic stabilization of systems \eqref{1.2} by means of an almost smooth feedback; (see \cite{art:3}, \cite{art:14} and \cite{art:16}). In order to provide the precise statement of \cite[Proposition 2]{art:19}, we first need to recall the following standard notations. For any pair of $C^{1} $ mappings $X: \mathbb{R}^{n} \to \mathbb{R}^{k}$, $Y: \mathbb{R}^{k} \to \mathbb{R}^{\ell } $ we adopt the notation $XY:=(DY)X$, $DY$ being the derivative of $Y$. By $[\cdot ,\cdot ]$ we denote the Lie bracket operator, namely, $[X,Y]=XY-YX$ for any pair of $C^{1} $ mappings $X,Y:\mathbb{R}^{ n} \to \mathbb{R}^{n} $. The precise statement of \cite[Proposition 2]{art:19} is the following. Assume that $f,g\in C^{2} $ and there exists a $C^{2} $, positive definite and proper function $V:{\mathbb R}^{n} \to {\mathbb R}^{+} $ such that the following implication holds:
\begin{equation} \label{1.3}
\begin{array}{l} {(gV)(x)=0,x\ne 0} \\ {\Rightarrow \left\{\begin{array}{c} {either\; (fV)(x)<0,({\rm ``Artstein-Sontag"\, condition})} \\ {or\; (fV)(x)=0;\; ([f,g]V)(x)\ne 0} \end{array}\right. } \end{array} 
\end{equation} 
Then system \eqref{1.2} is SDF-WGAS.

Proposition 2 of present work establishes that for systems \eqref{1.1} the same Lyapunov characterization of SDF-WGAS, originally proposed in \cite{art:17}, implies \textit{Semi-Global Asymptotic Stabilizability by means of a time-varying Sampled-Data Feedback }(SDF-SGAS), which is a stronger type of  SDF-WGAS. Proposition 3 is the main result of our present work. It constitutes a major generalization of \cite[Proposition 2]{art:19} mentioned above and provides a Lie algebraic sufficient condition for SDF-SGAS(WGAS) for the case of affine in the control systems \eqref{1.2}. This condition is much weaker than \eqref{1.3} and involves a particular Lie sub-algebra of the dynamics $f,g$ of the system \eqref{1.2}. 

The paper is organized as follows. Section II contains the  definitions of SDF-WGAS and SDF-SGAS and the statements of our results (Propositions 2 and 3). Section III contains the proofs of these results and in Section IV illustrative examples are provided. More results for 3-dimensional systems \eqref{1.2} are found in  \cite{art:arX}.

\section{Definitions and Main Results}
\label{S:2}
Consider system \eqref{1.1} and assume that $f:{\mathbb R}^{ n} \times {\mathbb R}^{ m} \to {\mathbb R}^{ n} $ is Lipschitz continuous. We denote by $x(\cdot )=x(\cdot ,s,x_{0} ,u)$ the trajectory of \eqref{1.1} with initial condition $x(s,s,x_{0} ,u)=x_{0} \in {\mathbb R}^{ n} $ corresponding to certain measurable and locally essentially bounded control $u:[s,T_{\max } )\to {\mathbb R}^{ m} $, where $T_{\max } =T_{\max } (s,x_{0} ,u)$ is the corresponding maximal existence time of the trajectory.
   
\begin{defn}
We say that system \eqref{1.1} is \textit{Weakly Globally Asymptotically Stabilizable by Sampled-Data Feedback (SDF-WGAS)}, if for every constant $\tau >0$ there exist mappings $T: {\mathbb R}^{ n} \backslash \{ 0\} \to  {\mathbb R}^{+} \backslash \{ 0\} $ satisfying 
\begin{equation} \label{2.1} T(x)\leq \tau ,\, \, \, \forall x\in  {\mathbb R}^{ n} \setminus \{ 0\}  \end{equation} 
and $k(t,x;x_{0} ): {\mathbb R}^{+} \times {\mathbb R}^{ n} \times {\mathbb R}^{ n} \to  {\mathbb R}^{ m} $ such that for any fixed $(x,x_{0} )\in {\mathbb R}^{n}\times \mathbb{R}^{n} $ the map $k(\cdot ,x;x_{0} ): {\mathbb R}^{+} \to  {\mathbb R}^{ m} $ is measurable and locally essentially bounded and such that for every $x_{0} \ne 0$ there exists a sequence of times
\begin{equation}\label{2.2} t_{1} :=0<t_{2} <t_{3} <\ldots <t_{\nu } <\ldots \; \, , \mathrm{with}\, \,  t_{\nu } \to \infty \end{equation}  
in such a way that the trajectory $x(\cdot )$ of the sampled-data closed loop system:
\begin{equation} \label{2.3} \begin{array}{l} {\dot{x}=f(x,k(t,x(t_{i} );x_{0} )),{\kern 1pt} {\kern 1pt} {\kern 1pt} {\kern 1pt} {\kern 1pt} {\kern 1pt} t\in [t_{i} ,\; t_{i+1} ),\; {\kern 1pt} {\kern 1pt} {\kern 1pt} {\kern 1pt} i=1,2,\ldots } \\ {\quad \quad \quad \quad \quad x(0)=x_{0} \in  {\mathbb R}^{ n} } \end{array} \end{equation} 
satisfies: 
\begin{equation} \label{2.4} t_{i+1} -t_{i} =T(x(t_{i} )),\; i=1,2,\ldots  \end{equation} 
and the following properties:
\begin{flalign}\text{Stability:}  &&
\begin{array}{c}\forall \varepsilon >0\Rightarrow \exists \delta =\delta (\varepsilon )>0: |x(0)|\le \delta\\
\Rightarrow |x(t)|\leq \varepsilon ,\; \forall t\geq 0\end{array} 
&&\label{2.5}
\end{flalign}
\begin{flalign}\text{Attractivity:}  && \mathop{\lim }\limits_{t\to \infty } x(t)=0,\; \forall x(0)\in  {\mathbb R}^{n}
&&\label{2.6}\end{flalign}
where ${\kern 1pt} \left|x\right|$ denotes the Euclidean norm of the vector $x$. 
\end{defn}

Next we give the Lyapunov characterization of SDF-WGAS proposed in \cite{art:17} and \cite{art:19}, that constitutes a generalization of the concept of the \textit{control Lyapunov function} (see Definition 5.7.1 in \cite{art:13}).

\begin{asmp}There exist a positive definite $C^{0} $ function $V: {\mathbb R}^{n} \to  {\mathbb R}^{+} $ and a function $a\in K$ (namely, $a(\cdot )$ is continuous, strictly increasing with $a(0)=0$) such that for every $\xi >0$ and $x_0\neq0$ there exists a constant $\varepsilon =\varepsilon(x_0)\in(0,\xi]$ and a measurable and locally essentially bounded control $u(\cdot,x_{0}) :[0,\varepsilon ]\to  {\mathbb R}^{m} $  satisfying
\begin{subequations} \label{2.7}
\begin{equation}\label{2.7a}
 V(x(\varepsilon ,0,x_{0} ,u(\cdot,x_{0}) ))<V(x_{0} );  
\end{equation}
\begin{equation}\label{2.7b}
 V(x(s,0,x_{0} ,u(\cdot,x_{0}) ))\le a(V(x_{0} )),\; \; \forall s\in [0,\varepsilon ]  
\end{equation}
\end{subequations}
\end{asmp}
The following result was established in \cite{art:17}.

\begin{prop}
Under Assumption 1, system \eqref{1.1} is SDF-WGAS.
\end{prop}

We now present the concept of  SDF-SGAS, which is a strong version of SDF-WGAS:

\begin{defn}We say that system \eqref{1.1} is \textit{Semi-Globally Asymptotically Stabilizable by Sampled-Data Feedback  (SDF-SGAS)}, if for every $R>0$ and for any given partition of times 
\begin{equation} \label{2.8} T_{1} :=0<T_{2} <T_{3} <\ldots <T_{\nu } <\ldots {\kern 1pt}  \, \, {\rm with}\, \, \,  T_{\nu } \to \infty  \end{equation} 
there exist a neighborhood $\Pi $ of zero with $B[0,R]:=\left\{x\in {\mathbb R}^{ n} :|x|\le R\right\}\subset \Pi $ and a map $k:{\mathbb R}^{+} \times \Pi \to  {\mathbb R}^{ m} $ such that for any $x\in \Pi $ the map $k(\cdot ,x): {\mathbb R}^{+} \to  {\mathbb R}^{ m} $ is measurable and locally essentially bounded and the trajectory $x(\cdot )$ of the sampled-data closed loop system
\begin{equation}  \begin{array}{c} {\dot{x}=f(x,k(t,x(T_{i} ))),\, \, t\in [T_{i} ,\; T_{i+1} ),\; \, \, i=1,2,\ldots } \\ { x(0)\in \Pi } \end{array} \label{2.9} \end{equation} 
satisfies:
\begin{flalign}\text{Stability:} &&
\begin{array}{c}\forall \varepsilon >0\Rightarrow \exists \delta =\delta (\varepsilon )>0:x(0)\in \Pi , \\
|x(0)|\le \delta\Rightarrow |x(t)|\leq \varepsilon ,\; \forall t\geq 0\end{array}
&&\label{2.10}
\end{flalign}
\begin{flalign}\text{Attractivity:}  && \mathop{\lim }\limits_{t\to \infty } x(t)=0,\; \forall x(0)\in \Pi&&\label{2.11}\end{flalign}
\end{defn}
\begin{rem}(i) It can be easily established that SDF-SGAS implies SDF-WGAS and the latter implies global asymptotic controllability at zero. 

(ii)  SDF-SGAS is stronger than the concept of sampled-data semi-global asymptotic stabilizability adopted in earlier works in the literature, because the selection of partition of the times in (2.8) is arbitrary. We also mention that, despite its semi-global nature, the advantage of SDF-SGAS is, according to Definition 2, that the partition of times in (2.8) and the corresponding control involved in (2.9) are independent of the initial state, while in Definition 1, the partition of times (2.2) and the corresponding control in (2.3) generally depend on the initial condition. This is an essential difference between SDF-SGAS and SDF-WGAS.
 
By exploiting the semi-global nature of Definition 2, particularly the requirement that \eqref{2.10} and \eqref{2.11} are valid for initial value $x(0)$ lying in a compact set, we can obtain the following proposition, which is one of the main results of the paper. Its proof is based on a generalization of the procedure employed in \cite{art:17} for the proof of Proposition 1.
\end{rem}

\begin{prop}\label{prp:2}
Under Assumption 1, system \eqref{1.1} is SDF-SGAS.
\end{prop}

We next present the precise statement of the central result of present work, which provides a Lie algebraic sufficient condition for SDF-SGAS(WGAS) for the affine in the control single-input system \eqref{1.2}. Assume that its dynamics $f$, $g$ are smooth ($C^{\infty}$) and let $ Lie\{ f,g\}$ be the Lie algebra generated by $\{f,g\}$. Let $L_{1}:=span\{f,g\}$ and $L_{i+1}:=span\{[X,Y],\, X\in L_{i},Y\in L_{1}\}$, $i=1,2,\ldots$ and for any nonzero $\Delta \in Lie\{f,g\}$ define:
\begin{equation} \label{2.12} order_{\left\{f,g\right\}} \Delta \left\{\begin{array}{l} {:=1,{\kern 1pt} {\rm if}\; \Delta \in L_{1} \backslash \{ 0\} } \\ \hspace{-0.52em}{\begin{array}{l} {:=k>1,{\rm if}{\kern 1pt} \Delta =\Delta _{1} +\Delta _{2} ,{\kern 1pt} \mathrm{with}\; \Delta _{1} \in L_{k} \backslash \{ 0\} } \\ {\quad \quad \quad \mathrm{and}{\kern 1pt} {\kern 1pt} \Delta _{2} \in span\{ \cup _{i=1}^{i=k-1} L_{i} \} } \end{array}} \end{array}\right.
 \end{equation} 

By exploiting the result of Proposition 2, the Campbell-Baker-Hausdorff (CBH) formula and applying a major extension of the proof of [18, Proposition 2] we get the following result for the case \eqref{1.2}, that constitutes the central result of present work.

\begin{prop}\label{prp:3}
For system \eqref{1.2} assume that there exists a smooth function $V: {\mathbb R}^{ n} \to {\mathbb R}^{+} $, being positive definite and proper, such that for every $x\neq 0$, either $(gV)(x)\neq 0$, or one of the following properties hold:\\ Either
\begin{equation}  (gV)(x )=0\Rightarrow (fV)(x )<0  \label{2.13}\end{equation}
or there exists an integer $N=N(x)\ge 1$ such that
\begin{subequations} \label{2.14}
\begin{equation}
\label{2.14a}
(gV)(x)=0,\; (f^{i} V)(x)=0,\; \; i=1,2,\ldots ,N\end{equation}
\begin{align}
(\Delta_{{1} }& \Delta_{{2} } \ldots \Delta_{{k} } V)(x)=0\nonumber \\
\forall \Delta _{{1} } , \Delta_{{2} }& ,\ldots , \Delta_{{k} } \in Lie\{ f,g\} \setminus \{ g\}\nonumber\\
\mathrm{with}&\, \, \sum\nolimits_{p=1}^{k}order_{\{f,g\}} \Delta_{{p} }  \le N \label{2.14b}
\end{align}
\end{subequations}
where $(f^iV)(x):=f(f^{i-1}V)(x)$, $i=2,3,\ldots$, $(f^1V)(x):=(fV)(x)$ and in such a way that one of the following properties hold:
\begin{flalign}\text{(P1)}  && (f^{N+1} V)(x)<0 && \label{2.15}\end{flalign}
(P2) $N$ is odd and
\begin{equation} \label{2.16} ([[\ldots [[f,\underbrace{g],g],\ldots ,g],g]}_{N\, \, \, \, }V)(x)\neq 0 \end{equation} 
(P3) $N$ is even and 
\begin{equation} \label{2.17} ([[\ldots [[f,\underbrace{g],g],\ldots ,g],g]}_{N\, \, \, \, }V)(x)<0 \end{equation} 
(P4) $N$ is an arbitrary positive integer with
\begin{subequations}\label{2.18}
\begin{align}(f^{N+1} V)(x)&=0,\;    \label{2.18a}\\
([[\ldots [[g,\underbrace{f],f],\ldots ,f],f]}_{N\, \, \, \, }&V)(x)\ne 0\; \label{2.18b}
 \end{align}
\end{subequations}
Then system \eqref{1.2} satisfies Assumption 1, hence, is SDF-SGAS and therefore SDF-WGAS.
\end{prop}

\begin{rem}(i) It should be pointed out, that the generalized concept of the control Lyapunov function given by Assumption 1, together with the result of Proposition 2, play a key role for the derivation of the Lie sufficient condition of Proposition 3; it should be emphasized here that the hypothesis of Proposition 3, guarantees the validity of Assumption 1 for system \eqref{1.2}, but it \textit{does not in general imply that $V$ involved in \eqref{2.13}-\eqref{2.18}  is a control Lyapunov function}, according to its standard definition in literature.

(ii) For the particular case of $N=1$, condition \eqref{2.14a} is equivalent to $(gV)(x)=0$ and $(fV)(x)=0$, the previous equality is equivalent to \eqref{2.14b} and obviously \eqref{2.16} is equivalent to $([f,g]V)(x)\ne 0$. It follows, according to the statement of Proposition 3, that, under \eqref{1.3}, system \eqref{1.2} is SDF-SGAS and thus SDF-WGAS; the latter conclusion, namely, that \eqref{1.3} implies SDF-WGAS, is the precise statement of \cite[Proposition 2]{art:19}. It turns out that Proposition 3 constitutes a generalization of the previously mentioned result in \cite{art:19}.

(iii) Statement of Proposition 3 is fulfilled under weaker regularity hypotheses for $f$, $g$ and $V$. Particularly, if we assume that $\bar{N}:=\sup\{N=N(x),\, x\neq 0\}<+\infty$, where $N=N(x)$ is the integer involved in (2.14)-(2.18), then the result of Proposition 3 holds under the assumption that $f$, $g$, $V\in C^{k}$ for certain integer $k>\bar{N}$. It also can be extended to multi-input affine in the control systems; for reasons of simplicity, only the single-input case is considered here.  
\end{rem}

\section{Proof of Main Results}
\begin{IEEEproof}[Proof of Proposition 2] 
Let $R$, $\rho $ be a pair of constants with $R>\rho \ge 0$ and define $S[\rho ,R):=\{x\in {\mathbb R}^{n} :\rho \le V(x)<R\}$. By exploiting \eqref{2.7a} and \eqref{2.7b} and applying similar arguments with those in proof of Proposition 1 in \cite{art:17}, it follows that for any $\xi >0$ there exists $\sigma \in (0,\xi ]$ such that for every $\varepsilon \in (0,\sigma ]$, a constant $L=L(\rho ,R)>0$ can be found in such a way that for every $t\ge 0$ and $x_{0} \in S[\rho ,R)$ there a exists a control $u(\cdot,x_0)$ (as determined in \eqref{2.7} with $\varepsilon$ as above) such that, if we define $u_{t}(s,x_{0}) :=u(s-t,x_{0})$, $s\in[t,t+\varepsilon ]$, the trajectory $x(\cdot ,\cdot ,x_{0} ,u_{t}(\cdot,x_0) )$ of \eqref{1.1} with $x(t,t,x_{0} ,u_{t}(\cdot,x_0) )=x_{0} $ satisfies:
\begin{subequations}\label{3.1}
\begin{align}
V(x(t+\varepsilon ,t,x_{0} ,u_{t}(\cdot,&x_0) ) ))\le V(x_{0} )-L;   \label{3.1a}\\
V(x(s,t,x_{0} ,u_{t}(\cdot,x_0) ))&\leq 2a(V(x_{0} )),\, \forall s\in \left[t,t+\varepsilon \right]\label{3.1b} 
\end{align}
\end{subequations}
Let $R>0$ arbitrary and let $\bar{R}>0$ be a constant such that $B[0,R]\subset S[0,\bar{R})$. Consider a partition of constants $\{ R_{n} ,\; n=1,2,\ldots \} $ with 
\begin{equation} \label{3.2}
R_{1} =\bar{R}, \, \, 0<R_{n+1} <R_{n} ,\, \, n=1,2,\ldots\, \,  \mathrm{with} \mathop{\lim }\limits_{n\to \infty } R_{n} =0\end{equation}
Also, let $\{ T_{\nu},\nu=1,2,\ldots \} $ be a  given partition of times satisfying \eqref{2.8}. For each $i=1,2,\ldots $ and constants $\varepsilon _{i} >0,i=1,2,...$ consider the following partition of times: 
\begin{equation}\label{3.3}P_{i} :=\left\{t_{i,1} :=0,t_{i,2} ,t_{i,3} ,\ldots \right\}\, \, \mathrm{with}\, \,  \mathop{\lim }\limits_{p\to \infty } t_{i,p} =\infty,\, \,  i=1,2,...\end{equation}
satisfying the following properties:  
\begin{subequations}
\label{3.4}
\begin{align}
0<t_{i,p} <&t_{i,p+1}   ;  \label{3.4a} \\
\{ T_{\nu} ,\nu=1,2,\ldots \} &\subset P_{i} \subset P_{i+1} ; \label{3.4b}\\
\varepsilon _{i} \ge t_{i,p+1} &-t_{i,p}>0 ,\; \; \forall i,p\in  {\mathbb N}   \label{3.4c}
\end{align}
\end{subequations}
Next we use the notation $u_{i,j}(\cdot,x_0):=u_{t_{i,j}}(\cdot,x_0)$. By using \eqref{3.1a} and \eqref{3.1b} with $\rho =R_{i+1} $ and $R=R_{i} $, $i=1,2,...$, we may find a constant  $L_{i} >0$, a partition of times and sufficiently small constant $\varepsilon _{i} >0$ such that \eqref{3.4} holds and simultaneously for $x_{0} \in S[R_{i+1} ,R_{i} )$ and for any pair of integers $(i,p)\in {\mathbb N}\times {\mathbb N}$, a control $u_{i,p}(\cdot;x_{0})  :[t_{i,p} ,t_{i,p} +\varepsilon _{i} ]\to  {\mathbb R}^{ m} $ can be found satisfying: 
\begin{subequations} \label{3.5}
\begin{align}
V(x(t_{i,p+1} ,&t_{i,p} ,x_{0} ,u_{i,p}(\cdot,x_{0}) ))\le V(x_{0} )-L_{i} ;   \label{3.5a}\\
V(x(s,t_{i,p} ,x_{0} ,&u_{i,p}(\cdot,x_{0}) ))\le 2a(V(x_{0} )),\forall s\in [t_{i,p} ,t_{i,p+1} ]   \label{3.5b}
\end{align}
\end{subequations}
The previous analysis asserts that, for given $\{ T_{\nu},\nu=1,2,\ldots \} $, a  partition of times \eqref{3.3} can be determined in such a way that \eqref{3.4a}, \eqref{3.4b} hold and simultaneously \eqref{3.5} is fulfilled, provided that $x_{0} \in S[R_{i+1} ,R_{i} )$. For each initial $x(0)\in \Pi :=S[0,R_{1} )$ consider the map $x(\cdot ):{\mathbb R}^{+} \to {\mathbb R}^{n} $ defined as follows:
\begin{subequations} \label{3.6}
\begin{equation} \label{3.6a}
\begin{array}{l} {\quad \; \; x(t)=x (t,t_{i,p} ,x(t_{i,p} ),u_{i,p}(\cdot,x_0) )} \\ {\forall t\in [t_{i,p} ,t_{i,p+1} ),\; \, x(t_{i,p} )\in S[R_{i+1} ,R_{i} ),\; \; i,p\in {\mathbb N}} \end{array}                             \end{equation}
where $x(\cdot,s,z,u)$ satisfies:
\begin{equation} \label{3.6b}
\dot{x }=f(x ,u),\, \, t\ge s, \; x (s,s,z,u)=z 
\end{equation}
\end{subequations}
 
An immediate consequence of  \eqref{3.3}, \eqref{3.4a}, \eqref{3.5} and \eqref{3.6} is the following fact:

\noindent \textbf{Fact 1:} The map $x(\cdot )$ as determined by \eqref{3.6} is well defined and satisfies:
\begin{subequations}\label{3.7}
\begin{align}
 V(x(t_{i,p+1} ))&\le V(x(t_{i,p} ))-L_{i} ; \label{3.7a}\\
V(x(s))\le 2a(V(x(t_{i,p} ))),&\forall s\in [t_{i,p} ,t_{i,p+1} ],\; i,p\in {\mathbb N} \label{3.7b}\\
\mathrm{provided\, \, that}\, \, & x(t_{i,p} )\in S[R_{i+1} ,R_{i} )\nonumber
\end{align}
\end{subequations}
and as a consequence of (3.7a) we get: 

\noindent \textbf{Fact 2:}
\begin{align} 
V(x(t_{k} ))\le V(x(t_{1}& ))-(k-1)\min \{ L_{j} ,{\kern 1pt} j=\nu ,\nu +1, \ldots \nonumber\\ \ldots,m\}, \, \,
\forall \, k,m,\nu &\in {\mathbb N}; m>\nu ,\, \, \, t_{i} \in P_{m} ,i=1,2,...,k:\nonumber\\ &t_{1} <t_{2} <\ldots <t_{k} \nonumber\\
\mathrm{provided\, that}\, \, x(&t_{1} ),x(t_{2} ),...,x(t_{k} )\in S[R_{m+1} ,R_{\nu } )\label{3.8}
\end{align}
and 
\begin{align}  V(x(t_{2} ))&\le V(x(t_{1} )), \forall t_{2} >t_{1} ;\; t_{2} ,t_{1} \in P_{\infty}:= \mathop{\bigcup }\limits_{i=1}^{\infty } P_{i} ,\nonumber\\
&x(t_{1} )\in \Pi  \label{3.9}\end{align} 
Moreover, by taking into account \eqref{3.4b}, \eqref{3.7b} and \eqref{3.9}, it follows:

\noindent\textbf{Fact 3}: For any $\tau  \in  P_{\infty} $ with $x(\tau )\in \Pi $, there exists a sequence $\{t_{k} ,k=1,2,...\}$ with $ t_{k} \in P_{\infty} $ and $t_{k+1} >t_{k} >\tau ,\; k=2,3,...,\; \; t_{1} :=\tau $ such that $\mathop{\lim }\limits_{k\to \infty } t_{k} =\infty $ and
\begin{equation} \label{3.10} V(x(s))\le 2a(V(x(t_{k} ))),\forall s\in [t_{k} ,t_{k+1} ) \end{equation} 
which by virtue of \eqref{3.9} implies:
\begin{equation} \label{3.11} V(x(s))\le 2a(V(x(t_{1} ))),\forall s\ge t_{1}  \end{equation} 
We next show that the map $x(\cdot )$ satisfies both \eqref{2.10} and \eqref{2.11}. Since $V$ is positive definite and proper, in order to establish \eqref{2.11}, it suffices to show that for initial nonzero $x(0)\in \Pi (=S[0,R_{1} ))$ and sufficiently small $\mu >0$ there exists a time $\tau \in  P_{\infty} $ such that 
\begin{equation} \label{3.12} V(x(t))\le \mu ,\forall t\ge \tau  \end{equation} 
Let $\theta ,\mu >0$ with $2a(\theta )<\mu $; $\theta \le R_{1} $ and let $m\in \mathbb{N}$ with 
\begin{equation} \label{3.13} R_{m+1} \le \theta <R_{m}  \end{equation} 
We claim that there exists $\bar{p}\in {\mathbb N}$ such that $t_{m,\bar{p}} \in P_{m} $ and
\begin{equation} \label{3.14} V(x(t_{m,\bar{p}} ))\le \theta  \end{equation} 
Indeed, otherwise we would have $\left\{x(t_{m,p} ):p=1,2,\ldots \right\}\cap S[0,R_{m+1} )=\emptyset $ and since $t_{m,p} \in P_{m} $, we obtain from \eqref{3.8} and \eqref{3.13} that
$R_{m+1}<V(x(t_{m,p} ))\le V(x(0))-(p-1)\min \{ L_{\nu } ,\, \nu=1,...,m\}$, $\forall p=1,2,\ldots
$,
a contradiction, hence, \eqref{3.14} is fulfilled. The latter, in conjunction with \eqref{3.10} and the definition of $\theta $ and $\mu $, implies $2a(V(x(t_{m,\bar{p}} )))\le 2a(\theta )<\mu $, which by virtue of \eqref{3.11}, asserts that for given $x(0)\in \Pi $ and sufficiently small constant $\mu >0$ there exists a time $\tau  \in  P_{\infty} $ such that  the map $x(\cdot )$ satisfies $V(x(t))\leq 2a(V(x(\tau )))<\mu $ for all $t\ge \tau $, which establishes \eqref{2.11}. Likewise, by using  \eqref{3.11} with $t_{1} =0$ we can establish that  \eqref{2.10} also holds for the map $x(\cdot )$. We are now in a position to establish that there exists a map $k: {\mathbb R}^{+} \times \Pi \to {\mathbb R}^{m} $ such that the trajectory of the sampled-data closed loop system \eqref{2.9} satisfies both \eqref{2.10} and \eqref{2.11}. Indeed, due to the first  inclusion of \eqref{3.4b}, for each  given $T_{i} $ and vector $z\in \Pi $ there exist times $t_{i_{k} ,p_{k} } \in P_{\infty} ,k=1,2,...,\nu $ and inputs $\omega _{k} :[t_{i_{k} ,p_{k} } ,t_{i_{k+1} ,p_{k+1} } )\to {\mathbb R}^{ m} $, $k=1,2,...,\nu -1$ such that  
\begin{subequations} \label{3.15}
\begin{equation} \label{3.15a}
\begin{array}{l}{t_{i_{k} ,p_{k} } <t_{i_{k+1} ,p_{k+1} } ;i_{k} \le i_{k+1} ;}\\ { i_{k} =i_{k+1} \Rightarrow p_{k+1} =p_{k} +1;} \\ { t_{i_{1} ,p_{1} } :=T_{i} ,\; \; t_{i_{\nu } ,p_{\nu } } :=T_{i+1} } \end{array} 
\end{equation}
\begin{equation} \label{3.15b}
\begin{array}{l} {x_{1} :=z;\, \; \omega _{1} (t):=u_{i_{1} ,p_{1} }(t,x_{1} ),t\in [t_{i_{1} ,p_{1} } ,t_{i_{2} ,p_{2} } ]} \\ {x_{2} :=x(t_{i_{2} ,p_{2} } ,t_{i_{1} ,p_{1} } ,x_{1} ,\omega _{1} );\; \; \omega _{2} (t):=u_{i_{2} ,p_{2}}(t,x_{2} )},\\ \quad \quad \quad \quad \quad t\in [t_{i_{2} ,p_{2} } ,t_{i_{3} ,p_{3} } ] \\ {x_{3} :=x(t_{i_{3} ,p_{3} } ,t_{i_{2} ,p_{2} } ,x_{2} ,\omega _{2} );\; \; \omega _{3} (t):=u_{i_{3} ,p_{3}}(t,x_{3}),}\\ \quad \quad \quad \quad \quad t\in [t_{i_{3} ,p_{3} } ,t_{i_{4} ,p_{4} } ] \\ {...} \\ {x_{\nu -1} :=x(t_{i_{\nu -1} ,p_{\nu -1} } ,t_{i_{\nu -2} ,p_{\nu -2} } ,x_{\nu -2} ,\omega _{\nu -2} );\; \;}\\ \quad \quad \omega _{\nu -1} (t):=u_{i_{\nu -1} ,p_{\nu -1}}(t,x_{\nu -1}),t\in [t_{i_{\nu -1} ,p_{\nu -1} } ,t_{i_{\nu } ,p_{\nu } } ] \end{array} 
\end{equation}
\end{subequations}
Then, if we define:
\begin{subequations} \label{3.16}
\begin{equation} \label{3.16a}
\begin{array}{l} {\phi _{i} (t,z):=\omega _{k} (t),t\in [t_{i_{k} ,p_{k} } ,t_{i_{k+1} ,p_{k+1} } ),\; z\in \Pi, } \\ {\quad \quad k=1,2,...,\nu -1, \quad t_{i_{1} ,p_{1} } =T_{i} ,\quad t_{i_{\nu } ,p_{\nu } } =T_{i+1} } \end{array} 
\end{equation}
\begin{equation} \label{3.16b}
k(t,z):=\phi _{i} (t,z),t\in [T_{i} ,T_{i+1} ),i=1,2,...,\; \; z\in \Pi   
\end{equation}
\end{subequations}
it is obvious that the map $x(\cdot )$ as defined in \eqref{3.6}  coincides  with the solution of the closed-loop \eqref{2.9} with $k:{\mathbb R}^{+} \times \Pi \to  {\mathbb R}^{ m} $ as defined by \eqref{3.15} and \eqref{3.16}, provided that their initial values at $t=0$ are the same. It turns out, according to stability analysis made for $x(\cdot )$, that \eqref{2.10} and \eqref{2.11} also hold  for the trajectory of the system \eqref{2.9} with $k: {\mathbb R}^{+} \times \Pi \to  {\mathbb R}^{ m} $ as defined above.
\end{IEEEproof}

\begin{IEEEproof}[Proof of Proposition 3]
Let $0\ne x_{0} \in {\mathbb R}^{ n} $ and suppose first that, either $(gV)(x_{0} )\ne 0$, or \eqref{2.13} is fulfilled with $x=x_0$, namely, $(gV)(x_{0} )=0$ and $(fV)(x_{0} )<0$. Then, in both cases above, there exists a constant input $u$ such that both \eqref{2.7a} and \eqref{2.7b} of Assumption 1 hold; particularly, for every sufficiently small $\varepsilon >0$ we have: 
\begin{equation} \label{3.17} V(x(s,0,x_{0} ,u))<V(x_{0} ),{\kern 1pt} \forall s\in (0,\varepsilon ] \end{equation} 
Assume next that there exists an integer $N=N(x_{0} )\geq 1$ satisfying \eqref{2.14}, as well as one of the properties (P1), (P2), (P3), (P4) with $x=x_0$. In order to derive the desired conclusion, we proceed as follows. Define: 
\begin{equation} \label{3.18} X:=f+u_{1} g,\, \, Y:=f+u_{2} g \end{equation} 
and for simplicity denote by $X_{t} (z)$ and $Y_{t} (z)$ the trajectories of the systems $\dot{x}=X(x)$ and $\dot{y}=Y(y)$, respectively, initiated at time $t=0$ from some $z\in  {\mathbb R}^{n} $. Also, for any constant $\rho>0$ define:
\begin{align} R(t):=(X_{\rho t} &\circ Y_{t} )(x_{0} ),t\ge 0,R(0)=x_{0} \label{3.19} \\
 m(t)&:=V(R(t)),t\ge 0  \label{3.20}
\end{align} 
and denote in the sequel by $\mathop{m}\limits^{(\nu )} (\cdot ),\, \nu =1,2,...$ its $\nu $-time derivative. We prove that, under previous assumptions concerning  the integer $N=N(x_{0} )$, there exist a constant $\rho=\rho(x_{0} )>0$ and  a pair of constant inputs $u_{1} $ and $u_{2} $ such that  $\mathop{m}\limits^{(n)}(0)=0,\,\, n=1,2,\ldots ,N$ and $\mathop{m}\limits^{(N+1)}(0)<0$. This would imply that $m(t)<m(0)=V(x_{0} )$ for every $t>0$ near zero and the latter in conjunction with \eqref{3.19} and \eqref{3.20} will lead to the validity of both inequalities \eqref{2.7a} and \eqref{2.7b} guaranteeing, according to Proposition 2, that \eqref{1.2} is SDF-SGAS. In order to get the desired result, we express the time derivatives $\mathop{m}\limits^{(\nu )}(0),\, \nu =1,2,...$ of the map $m(\cdot )$ in terms of the elements of the Lie algebra of $\left\{f,g\right\}$ and the function $V$ evaluated at $x_{0} $. We apply the CBH formula to the right hand side map of \eqref{3.19}. Then for  every $k\in {\mathbb N}$ we find:
\begin{align} 
\dot{R}(t) =&  \rho X(R(t))+(DX_{\rho t} Y)\circ X_{-\rho t} (R(t))\nonumber \\ 
=& \rho X(R(t))+Y(R(t))+\rho t[Y,X](R(t))\nonumber \\ 
&+\frac{\rho^{2} t^{2} }{2!} [[Y,X],X](R(t)) +\ldots \label{3.21}\\
&+\frac{\rho^{k} t^{k} }{k!} [...[[Y,\underbrace{X],X],\ldots ,X]}_{k}(R(t))+O(t^{k} )\nonumber 
\end{align} 
where $\mathop{\lim }\nolimits_{t\to 0^{+} } (O(t)/t)=0$. Let 
\begin{align} \label{3.22} 
A_{0}&:=  {\rho X+Y,}\nonumber \\
A_{\nu } &:= {[...[[Y,\underbrace{X],X],\ldots ,X]}_{\nu },\nu =1,2,...} 
\end{align} 
Notice that, since $A_{\nu } \in Lie\{ X,Y\} $, we may define, according to \eqref{2.12}, the order of each $A_{\nu } $ with respect to the Lie algebra of $\{ X,Y\} $; particularly, in our case, we have:
\begin{equation} \label{3.23} order_{_{\left\{X,Y\right\}} } A_{\nu } =\nu +1,\, \,  \forall \nu =0,1,2,\ldots  \end{equation} 
Now, \eqref{3.21} is rewritten: 
\begin{equation} \label{3.24} \dot{R}(t)=(A_{0} +\rho tA_{1} +\frac{1}{2!} \rho^{2} t^{2} A_{2} +\ldots +\frac{1}{k!} \rho^{k} t^{k} A_{k} )(R(t))+O(t^{k} ) \end{equation} 
thus, by invoking \eqref{3.20}, it follows that for any $k\in {\mathbb N}$ we have: 
\begin{align} \label{3.25} 
\mathop{m}\limits^{(1)} (t)=\sum_{i=0}^{k}\frac{\rho^it^i}{i!}(A_iV)(R(t))+O(t^{k} ) 
\end{align} 
Since we have assumed that $(fV)(x_{0} )=(gV)(x_{0} )=0$, it follows from \eqref{3.18}, \eqref{3.22} and \eqref{3.25} that 
\begin{equation} \label{3.26} \mathop{m}\limits^{(1)} (0)=0 \end{equation} 
From \eqref{3.24} and \eqref{3.25} we find:
\begin{align}
\mathop{m}\limits^{(2)} (t) =&\sum_{i=0}^{k}\frac{\rho^it^i}{i!} D(A_{i} V)(R(t))\dot{R}(t)+\sum_{i=1}^{k+1}\frac{\rho^it^{i-1}}{(i-1)!}(A_{i}V)(R(t))\nonumber\\
&+O(t^{k-1} )\nonumber\\
\in & (A_{0}^{2} V)(R(t)) +t\rho\, span\left\{A_{1} A_{0} V,A_{0} A_{1} V\right\}(R(t))\nonumber \\
 &+t^{2} \rho^{2} \, span\{A_{2} A_{0} V,A_{1}^{2} V,A_{0} A_{2} V\}(R(t)) \nonumber\\
 &+t^{3} \rho^{3} \, span\{A_{0} A_{3} V,A_{2} A_{1} V,A_{1} A_{2} V,A_{3} A_{0} V\}(R(t))\nonumber\\ 
 &+\ldots  +t^{k} \rho^{k} \, span\left\{A_{k} A_{0} V,A_{k-1} A_{1} V,..,A_{0} A_{k} V\right\}(R(t))\nonumber\\
 &+\rho(A_{1} V)(R(t))\nonumber\\
 &+span\{\rho^{2} tA_{2} V,\rho^{3} t^{2} A_{3} V,\ldots ,\rho^{k} t^{k-1} A_{k} V,\nonumber \\
 &\quad \quad  \rho^{k+1} t^{k} A_{k+1} V\}(R(t))+O(t^{k-1} )\label{3.27}
\end{align}

We show by induction that for every pair of integers $n,\, k$ with $2\le n\le k$, the n-time derivative $\mathop{m}\limits^{(n)} (\cdot )$ of $m(\cdot )$ satisfies:
\begin{align} &{\mathop{m}\limits^{(n)} (t)\in S_{n} (t,x_{0} ):=(A_{0}^{n} V)(R(t))}\nonumber\\
&{+\sum _{j=0}^{j=k} t^{j} span\left\{\begin{array}{l} {\rho^{r_{n}^{j} } (A_{i_{1}^{j} } A_{i_{2}^{j} } ...A_{i_{\nu }^{j} } V)(R(t)):\nu \ge 2;}\\{\sum _{s=1}^{\nu } order_{\left\{X,Y\right\}} A_{i_{s}^{j} } =n+j;\, } \\ {\, r_{n}^{j} =\sum _{s=1}^{\nu } i_{s}^{j} \in \left\{1,2,...,n+j-2\right\}} \end{array}\right\}}\nonumber\\
&+\rho^{n-1} (A_{n-1} V)(R(t))\nonumber \\ 
&+span\{ \rho^{n} t(A_{n} V)(R(t)),\rho^{n+1} t^{2} (A_{n+1} V)(R(t)),...,\nonumber\\
&\quad \quad \quad \quad \rho^{n+k-1} t^{k} (A_{k+n-1} V)(R(t))\} +O(t^{k-n+1} )\label{3.28} \end{align}
with $i_{1}^{j} ,i_{2}^{j} ,\ldots ,i_{\nu }^{j} \in  {\mathbb N}_{0} $, $j=0,1,2,\ldots ,k$. By taking into account \eqref{3.27}, it can be easily verified that inclusion \eqref{3.28} is indeed fulfilled for $n=2$. Suppose that \eqref{3.28} holds for some integer $n$, $2\le n<k$. We show that it is also fulfilled for $n=n+1\le k$. Indeed, from \eqref{3.28} the $(n+1)$-time derivative of $m(\cdot )$ is
\begin{align} 
&\mathop{m}\limits^{(n+1)} (t)=\frac{d}{dt} (\mathop{m}\limits^{(n)} (t))\in D(A_{0}^{n} V)(R(t))\dot{R}(t)\nonumber\\
&+\sum _{j=0}^{j=k} t^{j}  span\left\{\hspace{-0.44em}\begin{array}{l} {D(\rho^{r_{n}^{j} } A_{i_{1}^{j} } \ldots A_{i_{\nu }^{j} } V)(R(t)):\nu \ge 2;}\\ {\sum _{s=1}^{\nu } order_{\left\{X,Y\right\}} A_{i_{s}^{j} } =n+j;} \\ {r_{n}^{j} =\sum _{s=1}^{\nu } i_{s}^{j} \in \{ 1,2,\ldots ,n+j-2\} } \end{array}\hspace{-0.44em}\right\} \dot{R}(t)\nonumber\\
&{+\sum _{j=1}^{j=k} jt^{j-1}  span\left\{\begin{array}{l} {\rho^{r_{n}^{j} } (A_{i_{1}^{j} } \ldots A_{i_{\nu }^{j} } V)(R(t)):\nu \ge 2;}\\{\sum _{s=1}^{\nu } order_{\left\{X,Y\right\}} A_{i_{s}^{j} } =n+j;} \\ {r_{n}^{j}
 =\sum _{s=1}^{\nu } i_{s}^{j} \in \{ 1,2,\ldots ,n+j-2\} } \end{array}\hspace{-0.4em}\right\}} \nonumber\\
&+\rho^{n-1} D(A_{n-1} V)(R(t)) \dot{R}(t)\nonumber\\ 
&+span\{ \rho^{n} tD(A_{n} V)(R(t)),\rho^{n+1} t^{2} D(A_{n+1} V)(R(t)),...,\nonumber\\
&\quad \quad \rho^{n+k-1} t^{k} D(A_{k+n-1} V)(R(t))\}  \dot{R}(t)\nonumber\\
&+span\{ \rho^{n} (A_{n} V)(R(t)),\rho^{n+1} t(A_{n+1} V)(R(t)),...,\nonumber\\
&\quad \quad \rho^{n+j} t^{j} (A_{n+j} V)(R(t)),j=0,1,2,...,k\} +O(t^{k-n} )  \label{3.29}\end{align} 
Hence, by invoking \eqref{3.24} we have:
\begin{align} 
&\mathop{m}\limits^{(n+1)} (t) \in (A_{0}^{n+1} V)(R(t))\nonumber\\
&+span\left\{\rho^{q} t^{q} (A_{q} A_{0}^{n} V)(R(t)),q=1,...,n,n+1,...,k\right\} \nonumber\\
&+\hspace{-1.5em}\mathop{\sum \limits_{\begin{array}{l} \mbox{\small{{j=0,1,...,k}}} \\ \mbox{\small{{q=0,1,...,k}}}\\{j+q\le k} \end{array}} }\hspace{-1.52em} t^{j+q}  span\left\{\hspace{-0.4em}\begin{array}{l} {\rho^{r_{n}^{j} +q} (A_{q} A_{i_{1}^{j} } \ldots A_{i_{\nu }^{j} } V)(R(t)):\nu \ge 2;}\\{\sum _{s=1}^{\nu } order_{\left\{X,Y\right\}} A_{i_{s}^{j} } =n+j;} \\ {r_{n}^{j} =\sum _{s=1}^{\nu } i_{s}^{j} \in \{ 1,2,\ldots ,n+j-2\} } \end{array}\hspace{-0.4em}\right\}\nonumber\\
&+\sum _{j=1}^{j=k}t^{j-1}  span\left\{\hspace{-0.4em}\begin{array}{l} {j\rho^{r_{n}^{j} } (A_{i_{1}^{j} } \ldots A_{i_{\nu }^{j} } V)(R(t)):\nu \ge 2;}\\{\sum _{s=1}^{\nu } order_{\left\{X,Y\right\}} A_{i_{s}^{j} } =n+j;} \\ {r_{n}^{j} =\sum _{s=1}^{\nu } i_{s}^{j} \in \{ 1,2,\ldots ,n+j-2\} } \end{array}\hspace{-0.4em}\right\}\nonumber\end{align}\begin{align}
&+\rho^{n} (A_{n} V)(R(t))\nonumber \\
&+\rho^{n-1} span\{\rho^{q} t^{q} (A_{q} A_{n-1} V)(R(t));q=0,1,...,n,n+1,...,k\}\nonumber \\
&+span\{ \rho^{j+n-1+q} t^{j+q} (A_{q} A_{j+n-1} V)(R(t)),j=1,2,\ldots\nonumber\\
&\quad \quad \quad \ldots,n,n+1,...,k,q=0,1,...,k;j+q\le k\}\nonumber\\
&+span\{ \rho^{n+1} t(A_{n+1} V)(R(t)),...,\rho^{n+j} t^{j} (A_{n+j} V)(R(t)),\nonumber \\
&\quad \quad \quad j=1,2,...,k\} +O(t^{k-n} )\label{3.30} \end{align} 
Notice that each new term $t^{K} \rho^{L} A_{\tau _{1} } \ldots A_{\tau _{M} } V$ that appears above satisfies  
\begin{align}  &\sum _{s=1}^{s=M} order_{\left\{X,Y\right\}} A_{\tau _{s} } =(n+1)+K;\label{3.31}\\
 L&=\sum _{s=1}^{s=M} \tau _{s} \in \{ 1,2,\ldots ,(n+1)+K-2\}  \label{3.32}
 \end{align} 
For completeness, we note that for the terms $\rho^{q} t^{q} (A_{q} A_{0}^{n} V)$, $q=1,\ldots ,k$ it follows, by taking into account \eqref{3.28} and \eqref{3.29}, that $order_{\left\{X,Y\right\}} A_{q} +\sum _{s=1}^{s=n} order_{\left\{X,Y\right\}} A_{0} =(n+1)+q$ and obviously \eqref{3.32} holds as well. For the terms $t^{j+q} \rho^{r_{n}^{j} +q} (A_{q} A_{i_{1}^{j} } \ldots A_{i_{\nu }^{j} } V)$ we have: 
$order_{\left\{X,Y\right\}} A_{q} +\sum _{j=1}^{\nu } order_{\left\{X,Y\right\}} A_{i_{j}^{k} } =(n+1)+q+j$ 
and, since $r_{n}^{j} \in \{ 1,\ldots ,n+j-2\} $ as imposed in \eqref{3.30}, we have: 
$r_{n}^{j} +q\in \{ 1,2,\ldots ,n+q+j-2\} \subset  \{ 1,2,\ldots ,(n+1)+(q+j)-2\}$. 
Also, for  the terms $t^{j-1} \rho^{r_{n}^{j} } (A_{i_{1}^{j} } A_{i_{2}^{j} } \ldots A_{i_{\nu }^{j} } V)$ in \eqref{3.30} we have: 
$\sum _{j=1}^{\nu } order_{\left\{X,Y\right\}} A_{i_{j}^{k} } =(n+1)+j-1$
and obviously $r_{n}^{j} \in \{ 1,2,\ldots ,n+j-2\} \subset \{ 1,2,\ldots ,(n+1)+j-2\} $. Likewise, we handle the rest terms in the right hand side of \eqref{3.30} and show that both \eqref{3.31} and \eqref{3.32} hold. These conditions  imply  that the right hand set in \eqref{3.30} is included in $S_{n+1} (t,x_{0} )$, as the latter is defined in \eqref{3.28}, which guarantees that inclusion \eqref{3.28} holds for $n:=n+1$ and therefore is fulfilled for every pair of integers $k\ge n\ge 2$. It follows from \eqref{3.27} and \eqref{3.28} that
\begin{equation} \label{3.33} \mathop{m}\limits^{(2)} (0)=(A_{0}^{2} V)(x_{0} )+(\rho A_{1} V)(x_{0} ) \end{equation} 
for the case $n=2$ and generally for $n\ge 2$:
\begin{align}&\mathop{m}\limits^{(n)} (0) \in (A_{0}^{n} V)(x_{0} )\nonumber\\
&+span\left\{\hspace{-0.3em}\begin{array}{l} {\rho^{r_{n}^{0} } (A_{i_{1}^{0} } A_{i_{2}^{0} } ...A_{i_{\nu }^{0} } V)(x_{0} ):{\kern 1pt} \nu \ge 2;}\\{i_{1}^{0} ,i_{2}^{0} ,...i_{\nu }^{0} \in {\mathbb N}_{0} ;\sum _{j=1}^{\nu } order_{\left\{X,Y\right\}} A_{i_{j}^{0} } =n;} \\ {r_{n}^{0} =\sum _{j=1}^{\nu } i_{j}^{0} \in \left\{1,2,...,n-2\right\}} \end{array}\right\}\nonumber\\
&+\rho^{n-1} (A_{n-1} V)(x_{0} )\label{3.34} \end{align} 
By taking into account definition \eqref{3.18} of the vector fields $X$ and $Y$ and by setting
\begin{equation} \label{3.35} u_{2} =-\rho u_{1}, \, \rho>0  \end{equation} 
we get
\begin{align} 
A_{0} &=(\rho+1)f, \, \,  A_{1} =(\rho+1)u_{1} [f,g]\nonumber\\
A_{2} &=  (\rho+1)(u_{1}^{2} [[f,g],g]-u_{1} [[g,f],f])\nonumber \\[-0.5em] 
&\vdots \nonumber\\[-0.5em] 
A_{n} &= (\rho+1)u_{1}^{n} [\ldots [[f,\underbrace{g],g],\ldots ,g]}_{n{\kern 1pt} {\kern 1pt}}\nonumber\\
&+(\rho+1)u_{1}^{n-1} ([[[\ldots [f,\underbrace{g],\ldots ,g],g}_{n-1{\kern 1pt} {\kern 1pt} }],f]\nonumber \\
  & +[[[\ldots [f,\underbrace{g],\ldots ,g}_{n-2{\kern 1pt} {\kern 1pt} }],f],g]+\ldots +[\ldots [[[f,g],f],\underbrace{g]\ldots ,g]}_{n-2{\kern 1pt} {\kern 1pt} })\nonumber\\
 &+\ldots +(\rho+1)u_{1}^{2} ([[[\ldots [[f,g],\underbrace{f],\ldots ,f],f}_{n-2{\kern 1pt} {\kern 1pt} }],g]\nonumber\end{align}\begin{align}
 &+[[[...[[f,g],\underbrace{f],... ,f}_{n-3{\kern 1pt} {\kern 1pt} }],g],f]
 +... +[[... [[[f,g],g],\underbrace{f]... ,f],f}_{n-2{\kern 1pt} {\kern 1pt} }])\nonumber \\
 &-(\rho+1)u_{1} [... [[g,\underbrace{f],f],... ,f}_{n{\kern 1pt} {\kern 1pt} }],\; n=3,4,...  \label{3.36}
\end{align}
Obviously, \eqref{3.36} implies:
\begin{align}  
A_{k} \in span\{\Delta \in Lie&\left\{f,g\right\}\setminus \{g\}:order_{\{f,g\}}\Delta  =k+1\}\nonumber \\
&k=0,1,2,\ldots\label{3.37}
\end{align} 
Also, we recall from \eqref{3.23} and \eqref{3.34} that 
$r_{n}^{0} =\sum _{s=1}^{\nu } i_{s}^{0} \in \{ 1,2,\ldots ,n-2\} $ and $\sum _{j=1}^{\nu } order_{\{X,Y\}} A_{i_{j}^{0}}=r_{n}^{0}+\nu =n$ with $\nu \ge 2$
and  therefore $\nu \le n-1$. By \eqref{3.34}-\eqref{3.37} and the previous facts we get:
\begin{align} \mathop{m}\limits^{(n)} (0)&\in (\rho+1)^{n} (f^{n} V)(x_{0} )+u_{1}^{} \pi _{1} (\rho,\rho+1;x_{0} )\nonumber \\ &+span\left\{u_{1}^{k} \pi _{k} (\rho,\rho+1;x_{0} ),k=2,...,n-2\right\}\nonumber \\  
&+\rho^{n-1} (\rho+1)u_{1}^{n-1} ([\ldots [[f,\underbrace{g],g],\ldots ,g}_{n-1\, \, {\kern 1pt} {\kern 1pt} }]V)(x_{0} ) \nonumber \\ 
&-\rho^{n-1} (\rho+1)u_{1} ([\ldots [[g,\underbrace{f],f],\ldots ,f}_{n{\kern 1pt} -1\, \, \, \, {\kern 1pt} }]V)(x_{0} ) \label{3.38} 
\end{align} 
for $n=2,3,...$ and for certain smooth functions $\pi _{k} : {\mathbb R}^{2} \times {\mathbb R}^{n} \to  {\mathbb R},\,\,  k=1,2,\ldots ,n-2$ satisfying the following properties:

\noindent(S1) For every $x_{0} \in  {\mathbb R}^{n} $, each map $\pi _{k} (\alpha ,\beta ;x_{0} ): {\mathbb R}^{2} \to  {\mathbb R}$ is a polynomial with respect to the first two variables  in such a way that
\begin{equation} \label{3.39} \begin{array}{l} { span\{ \pi _{k} (\alpha ,\beta ;x_{0} ),\, {\kern 1pt} k=1,2,\ldots ,n-2\} \subset } \\ {span\{ (\Delta _{{1} }  \Delta _{{2} } ...\Delta _{{i} } V)(x_{0} );\; {i} \in {\mathbb N} ,\; } \\ {\quad {\kern 1pt} \; \; \; \Delta _{{1} }  ,\Delta _{{2} } ,...,\Delta _{{i} } \in Lie\{ f,g\} \backslash \{ g\} ;{ \sum _{j=1}^{j=i} order_{\left\{f,g\right\}} \Delta _{{j} } =n\, \, \} }} \end{array} \end{equation} 
(S2) For each $x_{0} \in  {\mathbb R}^{n} $ there exist integers $\lambda _{i}$, $\mu _{i}$, $i=1,2,...,\; L\in {\mathbb N}$ with $1\le \lambda _{i} \le n-2$, $2\le \mu _{i} \le n-1$ such that the map $\pi _{1} (\alpha ,\beta ;x_{0} ): {\mathbb R}^{2} \to \mathbb{R}$ satisfies
$ \pi _{1} (\alpha ,\beta ;x_{0} )\in span\left\{\alpha^{\lambda _{1} } \beta ^{\mu _{1} } ,\alpha^{\lambda _{2} } \beta ^{\mu _{2} } ,...,\alpha^{\lambda _{L} } \beta ^{\mu _{L} } \right\}$.
The latter implies that for each fixed $x_{0} \in  {\mathbb R}^{n} $ the polynomials $\pi _{1} (\rho,\rho+1;x_{0} )$ and $-\rho^{n-1} (\rho+1)([\ldots [[g,\underbrace{f],f],\ldots ,f}_{n -1\, \, \, \,  }]V)(x_{0} )$ are linearly independent, provided that
\begin{equation} \label{3.41} ([[\ldots [[g,\underbrace{f],f],\ldots ,f],f]}_{n-1\, \, \, \, }V)(x_{0} )\ne 0\;  \end{equation} 
If we define:
\begin{align} \label{3.42} \xi _{n} (\rho;x_0):=&\pi _{1} (\rho ,\rho+1;x_{0} )\\
&-\rho^{n-1} (\rho+1)([\ldots [[g,\underbrace{f],f],\ldots ,f}_{n{\kern 1pt} -1\, \, \, \, {\kern 1pt} }]V)(x_{0} )\nonumber \end{align} 
the inclusion \eqref{3.38} is rewritten:    
\begin{align} \mathop{m}\limits^{(n)} (0)&\in (\rho+1)^{n} (f^{n} V)(x_{0} )+u_{1}^{} \xi _{n} (\rho;x_{0} )\nonumber \\
&+span\left\{u_{1}^{k} \pi _{k} (\rho,\rho+1;x_{0} ),k=2,...,n-2\right\}\nonumber \\ 
&+\rho^{n-1} (\rho+1)u_{1}^{n-1} ([\ldots [[f,\underbrace{g],g],\ldots ,g}_{n-1\, \, {\kern 1pt} {\kern 1pt} }]V)(x_{0} ) \label{3.43} \end{align} 
and a  constant $\rho=\rho(x_{0} )>0$ can be found with
\begin{equation} \label{3.44} \begin{array}{l} {\xi _{n} (\rho;x_{0} )\ne 0} \end{array} \end{equation} 
provided that \eqref{3.41} holds. Suppose now that there exists an integer $N=N(x_{0} )\ge 1$ satisfying  \eqref{2.14}, as well as one of the properties  (P1), (P2), (P3), (P4) with $x=x_0$. By \eqref{3.26} and by taking into \nopagebreak[0]account \eqref{2.14}, \eqref{3.38} and \eqref{3.39} it  follows:
\begin{equation} \label{3.45} \mathop{m}\limits^{(n)} (0)=0,{\kern 1pt} {\kern 1pt} n=1,2,\ldots ,N \end{equation} 
and we distinguish four cases:

\noindent\textbf{Case 1:} \textit{\eqref{2.15} holds with $x=x_0$}. Then by using \eqref{3.43} with $n:=N+1$ and by setting  $u_{1} =0$ we find that for all $\rho>0$ it holds:
\begin{equation} \label{3.46} \mathop{m}\limits^{(N+1)} (0)<0 \end{equation} 

\noindent\textbf{Case 2: }\textit{$N$ is odd and \eqref{2.16} holds with $x=x_0$. }We again invoke \eqref{3.43} with $n:=N+1$ and our assumption that $N$ is odd. It follows that for every $\rho>0$ there exists a constant  $u_{1} =u_{1} (x_{0} )$, with $|u_1|$ sufficiently large, such that again \eqref{3.46} is fulfilled.

\noindent\textbf{Case 3:} \textit{$N$ is even and \eqref{2.17} holds with $x=x_0$}. Then, as in the previous case, by using \eqref{3.43} with $n:=N+1$ it follows that, for any choice of $\rho>0$ and for any  sufficiently large constant $u_{1} =u_{1} (x_{0} )>0 $, the desired \eqref{3.46} holds.

\noindent\textbf{Case 4:} \textit{$N$ is arbitrary and both \eqref{2.18a}  and \eqref{2.18b} are satisfied with $x=x_0$}. Then, due to assumption \eqref{2.18b}, it follows that \eqref{3.41} is fulfilled with $n:=N+1$, therefore there exists a constant $\rho=\rho(x_{0} )>0$ satisfying \eqref{3.44} with $n:=N+1$. By invoking again \eqref{3.43} with $n:=N+1$ and by taking into account assumption \eqref{2.18a}, it follows that for this $\rho$ above there exists a sufficiently small constant $u_{1} =u_{1} (x_{0} )\ne 0$ such that \eqref{3.46} holds.

It follows, by taking into account \eqref{3.19}, \eqref{3.20}, \eqref{3.35}, \eqref{3.45} and \eqref{3.46}, that in all previous cases, there exists a constant $u_{1} $ such that, if for any $t>0$  we define: 
\begin{equation} \label{3.47} \omega(s;t,x_{0}):=\left\{\begin{array}{c} {u_{2} =-\rho u_{1} ,{\kern 1pt} s\in [0,t]} \\ {u_{1} ,{\kern 1pt} s\in (t,t+\rho t]} \end{array}\right.  \end{equation} 
 with $\rho=\rho(x_0):=1$ for the Cases 1, 2 and 3 and $\rho=\rho(x_0)$ as considered in the Case 4, then for every sufficiently small $\sigma=\sigma(x_0) >0$ we have $ m(t)<m(0)$, $\forall t\in (0,\sigma ]$, where $m(t):=V((X_{\rho t} \circ Y_{t} )(x_{0} ))=V(x(t+\rho t,0,x_{0} ,\omega(\cdot;t,x_{0} ))$ and $x(\cdot ,0,x_{0} ,\omega(\cdot;t,x_{0} ) )$ is the trajectory of \eqref{1.2} corresponding to the input $\omega(\cdot;t,x_{0} ) $. Equivalently:
\begin{equation} \label{3.48} V(x(t,0,x_{0} ,\omega(\cdot;\tfrac{t}{1+\rho},x_{0} ) ))<V(x_{0} ){\kern 1pt} ,\forall t\in(0,\tfrac{\sigma}{1+\rho}] \end{equation} 
Since the constant $\rho=\rho(x_{0} )$ is independent of $t$, we may pick $\varepsilon \in (0,\sigma]$ sufficiently small in such a way that, if we define 
$u(\cdot,x_0):=\omega(\cdot;\tfrac{\varepsilon}{1+\rho},x_0)$, 
inequality in \eqref{3.48} holds for $t:=\varepsilon $, namely, $V(x(\varepsilon,0,x_{0},$ $u(\cdot,x_{0})))<V(x_{0})$ and simultaneously $ V(x(s,0,x_{0} ,u(\cdot,x_{0}) ))$ $\le 2V(x_{0} )$, $ \forall s\in (0,\varepsilon ] $.
We conclude, by taking into account \eqref{3.17} and previous inequalities, that for every $x_{0} \ne 0$ and $\xi>0$, there exist $\varepsilon=\varepsilon(x_0)\in(0,\xi]$ and a measurable and locally essentially bounded control $u(\cdot,x_{0}) :[0,\varepsilon ]\to {\mathbb R}$ such that \eqref{2.7a} and \eqref{2.7b} of Assumption 1 hold with $a(s):=2s$. Therefore, according to Proposition 2, \eqref{1.2} is SDF-SGAS.   
\end{IEEEproof}

\section{Examples}\label{S:3}
The following examples illustrate the nature of Proposition 3. The first example below generalizes Example 2 in \cite{art:17}. 
\begin{exmp}
For the planar case: 
$ \dot{x}_{1} =F(x_{1} ,x_{2} ),\dot{x}_{2} =u,\, \, (x_{1} ,x_{2} )\in {\mathbb R}^{2}$,
where $F:{\mathbb R}^{2} \to {\mathbb R}^{+} $ is $C^{\infty } $, assume that for every $x_{1} \ne 0$, either $x_{1} F(x_{1} ,0)<0 $, or there exists an integer $N=N(x_1)\ge 1$ with $ \frac{\partial ^{i} F}{\partial x_{2}^{i} } (x_{1} ,0)=0$, $i=0,1,...,N-1 $ such that one of the following properties hold: (H1) $N$ is odd and $\frac{\partial ^{N} F}{\partial x_{2}^{N} } (x_{1} ,0)\ne 0$; (H2) $N$ is even and $x_{1} \frac{\partial ^{N} F}{\partial x_{2}^{N} } (x_{1} ,0)<0$. Then by setting $x:=(x_{1}^{} ,x_{2}^{} )^{T} $, $V(x):={\tfrac{1}{2}} (x_{1}^{2} +x_{2}^{2} )$, $f(x):=(F(x_{1} ,x_{2} ),0)^{T} $ and $g(x):=(0,1)^{T} $ it follows that for those $x\neq 0$ for which $(gV)(x)=0$, either \eqref{2.13} holds, or \eqref{2.14} together with one of the properties (P2), (P3) of Proposition 3 are fulfilled, hence, the system is SDF-SGAS.
\end{exmp}
\begin{exmp}
Consider the system $\dot{x}_{1} =x_{2} \alpha(x_{3} ),\,\,\dot{x}_{2}=-x_{1}\beta(x_{3} ),\,\,\dot{x}_{3}=u,\, \, (x_{1} ,x_{2},x_3 )\in {\mathbb R}^{3}$, where $\alpha(\cdot ),\beta(\cdot )\in C^{2} ( {\mathbb R}, {\mathbb R})$, which satisfy $\alpha(0)=\beta(0)\ne 0$ and $\mathop{\alpha}\limits^{(1)} (0)\ne \mathop{\beta}\limits^{(1)} (0)$, where $\mathop{\alpha}\limits^{(1)} (\cdot )$ and $\mathop{\beta}\limits^{(1)} (\cdot )$ denote the first derivatives of the functions $\alpha(\cdot )$ and $\beta(\cdot )$, respectively.  Define $x:=(x_{1}^{} ,x_{2}^{} ,x_{3} )^{T} $, $f(x):=(x_{2} \alpha(x_{3} ),-x_{1} \beta(x_{3} ),0)^{T} $, $g(x):=(0,0,1)^{T} $ and $V(x):={\tfrac{1}{2}} (x_{1}^{2} +x_{2}^{2} +x_{3}^{2} )$. Let $x\neq0$ and suppose that $(gV)(x)=x_3=0$. It follows that $(fV)(x)=(f^2V)(x)=(f^3V)(x)=0$. We distinguish two cases. The first is $([f,g]V)(x)\neq 0$, $x\neq 0$, which in conjunction with the previous equalities, assert that  \eqref{2.14} and \eqref{2.16} of Proposition 3 hold with $N=1$. The second case is $([f,g]V)(x)=0$, which, in conjunction with $(gV)(x)=x_3=0,\,x\neq0$ and hypotheses imposed for the terms $\alpha(\cdot)$ and $\beta(\cdot)$, guarantees that $([[g,f],f]V)(x)\neq0$, namely \eqref{2.18b} holds with $N=2$. It is also obvious that in this case, condition \eqref{2.14} is fulfilled as well with $N=2$. It turns out, according to the statement of Proposition 3, that the system is SDF-SGAS. 
\end{exmp}


\begin{thebibliography}{30}
\bibitem{art:1}{F. Ancona and A. Bressan, ``Patchy vector fields and asymptotic stabilization'', \textit{ESAIM-COCV}, vol. 4, pp. 445-471, 1999.}

\bibitem{art:2}{A. Anta and P. Tabuada, ``To sample or not to sample: self-triggered control for nonlinear systems'', \textit{IEEE Trans. Autom. Control}, vol. 55, pp. 2030-2042, 2010.}

\bibitem{art:3}{Z. Artstein, ``Stabilization with relaxed controls'', \textit{Nonlinear Analysis TMA}, vol. 7, pp. 1163-1173, 1983.}

\bibitem{art:4}{A. Bacciotti and L. Mazzi, ``From Artstein-Sontag Theorem to the min-projection strategy, \textit{Trans. of the Institute of Measurement and Control}, vol. 32, no. 6, pp. 571-581, 2010.}

\bibitem{art:6}{F.H. Clarke, Y.S. Ledyaev, E.D. Sontag and A.I. Subbotin, ``Asymptotic controllability implies feedback stabilization'', \textit{IEEE Trans. Autom. Control}, vol. 42, no. 10, pp. 1394-1407, 1997.}

\bibitem{art:7}{L. Gr\"{u}ne and D. Ne\v{s}i\'{c}, ``Optimization based stabilization of sampled-data nonlinear systems via their approximate discrete-time models'', \textit{SIAM J. Control  Optim.}, vol. 42 , pp. 98-122, 2003. }

\bibitem{art7.1}{I. Karafyllis, ``Stabilization by means of time-varying hybrid feedback'', \textit{MCSS}, vol. 18, pp. 236-259, 2006}

\bibitem{art:8}{N. Marchand and M. Alamir, ``Asymptotic controllability implies continuous discrete-time feedback stabilization'', In: Nonlinear Control in the Year 2000, vol. 2, Springer, Berlin, Heidelberg, New York, 2000.}

\bibitem{art:8.1}{H. Michalska and M.Torres-Torriti, ``A geometric approach to feedback stabilization of nonlinear systems with drift'', \textit{Systems and Control Lett.}, vol. 50, no. 4, pp. 303-318, 2003.}

\bibitem{art:10}{D. Ne\v{s}i\'{c} and A.R. Teel, ``A framework for stabilization of nonlinear sampled-data systems based on their approximate discrete-time models'', \textit{IEEE Trans. Autom. Control}, vol. 49, pp. 1103-1122, 2004.}

\bibitem{art:11}{C. Prieur, ``Asymptotic controllability and robust asymptotic stabilizability'', \textit{SIAM J. Control Optim.}, vol. 43, pp. 1888-1912, 2005. }

\bibitem{art:11b}{C. Prieur, R. Goebel and A.R. Teel, ``Hybrid feedback control and robust stabilization of nonlinear systems'', \textit{SIAM J. Control Optim.}, vol. 43, pp. 1888-1912, 2005.}

\bibitem{art:12}{H. Shim and A.R. Teel, ``Asymptotic controllability and observability imply semiglobal practical asymptotic stabilizability by sampled-data output feedback'', \textit{Automatica}, vol. 39, pp. 441-454, 2003.}

\bibitem{art:13}{E.D. Sontag, \textit{Mathematical control theory}, 2$^{nd}$ edn., Springer, Berlin, Heidelberg, New York, 1998.}

\bibitem{art:14}{E.D. Sontag, ``A `universal' construction of Artstein’s theorem on nonlinear stabilization'', \textit{Systems and Control Lett.}, vol. 13, pp. 117-123, 1989.}

\bibitem{art:16}{J. Tsinias, ``Sufficient Lyapunov-like conditions for stabilization'', \textit{Math. Contr. Sign. Syst.}, vol. 2, pp. 343-357, 1989.}

\bibitem{art:17}{J. Tsinias, ``Remarks on asymptotic controllability and sampled-data feedback stabilization for autonomous systems'', \textit{IEEE Trans. Autom. Control}, vol. 55, pp. 721-726, 2010.}

\bibitem{art:19}{J. Tsinias, ``New results on sampled-data feedback  stabilization for autonomous nonlinear systems'', \textit{Systems and Control Lett.}, vol. 61, pp. 1032-1040, 2012.    }
\bibitem{art:arX}{J. Tsinias and D. Theodosis, ``Sufficient Lie Algebraic Conditions for Sampled-Data Feedback Stabilization'', arXiv:1506.00078 [math.OC], to appear for presentation in \textit{54th IEEE CDC}, 2015.}
\end{thebibliography}
\end{document}